\documentclass[draftcls,onecolumn]{IEEEtran}
\ifCLASSINFOpdf
\usepackage[pdftex]{graphicx}
\else
\fi
\usepackage{amsmath}

\usepackage{algorithmic}
\ifCLASSOPTIONcompsoc
 \usepackage[caption=false,font=normalsize,labelfont=sf,textfont=sf]{subfig}
\else
 \usepackage[caption=false,font=footnotesize]{subfig}
\fi
\usepackage{url}

\usepackage{amssymb}
\usepackage{verbatim}

\newtheorem{theorem}{Theorem}
\newtheorem{lemma}{Lemma}

\newtheorem{proposition}{Proposition}

\DeclareMathOperator*{\argmin}{arg\,min}

\DeclareMathOperator{\E}{\mathbb{E}}
\newcommand*{\QED}{\hfill\IEEEQEDopen}%
\newcommand{\ee}{\mathrm{e}}
\renewcommand{\Re}{\mathbb{R}}

\newcommand{\cc}{\mathbf{c}}
\newcommand{\gt}{\mathbf{g}_t}
\newcommand{\mub}{\boldsymbol{\mu}}

\newcommand{\vb}{\mathbf{v}}
\newcommand{\Tr}{^{\mathrm{T}}}
\newcommand{\norm}[1]{\left\lVert#1\right\rVert}
\newcommand{\Sone}{\mathbb{S}_1}
\newcommand{\BR}{B_\mathcal{R}}
\newcommand{\KK}{\mathcal{K}}
\newcommand{\RR}{\mathcal{R}}

\newcommand{\xx}{\mathbf{x}}

\begin{document}
%
\title{Setpoint Tracking with Partially Observed Loads}
%
%
%

\author{Antoine~Lesage-Landry,~\IEEEmembership{Student Member,~IEEE,}
        and~Joshua~A.~Taylor,~\IEEEmembership{Member,~IEEE}
\thanks{This work was funded by the Fonds de recherche du Qu\'ebec -- Nature et technologies and the Natural Sciences and Engineering Research Council of Canada.}
\thanks{A. Lesage-Landry and J.A. Taylor are with The Edward S. Rogers Sr. Department of Electrical \& Computer Engineering, University of Toronto, Toronto,
Ontario, Canada, M5S 3G4. e-mail: \{\texttt{alandry@ece.,josh.taylor@}\}\texttt{utoronto.ca.}}}

%
%

\markboth{}{}%

%



\maketitle

\begin{abstract} 
We use online convex optimization (OCO) for setpoint tracking with uncertain, flexible loads. We consider full feedback from the loads, bandit feedback, and two intermediate types of feedback: partial bandit where a subset of the loads are individually observed and the rest are observed in aggregate, and Bernoulli feedback where in each round the aggregator receives either full or bandit feedback according to a known probability. We give sublinear regret bounds in all cases. We numerically evaluate our algorithms on examples with thermostatically controlled loads and electric vehicles.
\end{abstract}


%
\IEEEpeerreviewmaketitle

\section{Introduction}
\IEEEPARstart{D}{emand} response (DR) is an important source of flexibility for the electric power grid~\cite{callaway2011achieving,siano2014demand,taylor2016power}. In this work, we use online convex optimization (OCO) to design algorithms for tracking setpoints with uncertain flexible loads in DR programs.

Uncertainty is a key challenge in DR~\cite{taylor2015uncertainty}. This uncertainty arises from weather, human behavior and unknown load models. In addition, loads are time-varying, e.g., a heater may not have much flexibility during cold nights. The uncertain, time-varying nature of loads means that load aggregators must deploy loads for DR to observe their capabilities.

We base our setpoint tracking model on OCO~\cite{shalev2011online,hazan2016introduction}. This allows us dispatch loads for DR without precise knowledge of their responses by relying only on information from previous rounds. We invoke theoretical bounds to guarantee the performance of the approach.

We focus on controllable loads that can both increase and decrease their power consumption when instructed. In each round, the load aggregator sends adjustment signals to loads so that the aggregate power consumption tracks a setpoint. We use a sparsity regularizer and a mean regularizer to reduce the number of dispatched loads and the impact on the loads,  respectively. To minimize the setpoint tracking objective and regularizers, Composite Objective MIrror Descent (\texttt{COMID})~\cite{duchi2010composite}, an OCO algorithm, is used to compute the adjustment signal in each round.

\subsection{Related work}
Online learning~\cite{shalev2011online,bubeck2012regret} and online convex optimization~\cite{zinkevich2003online,bubeck2011introduction,hazan2016introduction} have already seen a number of applications in DR. Several variants of the multi-armed bandit framework have been used to curtail flexible loads in~\cite{taylor2014index,wang2014adaptive,kalathil2015online,bandyopadhyay2016planning}. Reference~\cite{lesagelandry2017learning} used adversarial bandits to shift load while learning load parameters. 

Online learning has also been used in models of price-based DR. Reference~\cite{jia2013retail} uses a continuum-armed bandit to do real-time pricing of price responsive dynamic loads. Reference~\cite{soltani2015real} uses OCO and conditional random fields to predict the price sensitivity of Electric Vehicles (EVs). This model was then used as input to compute real-time prices. Using OCO,~\cite{kimonline} develops real-time pricing algorithms to flatten the aggregate load consumption. They later apply their algorithm to EVs charging. OCO was also used in~\cite{ma2016distributed} to flatten the aggregated power consumption using EVs charging scheduling. Ledva et al.~\cite{ledva2015inferring} used OCO to identify the controllable portion of demand in real-time from aggregate measurements.

\subsection{Contributions}
Reference~\cite{kimonline} is the most closely related to ours. In~\cite{kimonline}, \texttt{COMID} is used to set prices with the objective of flattening load. We differ in that we use direct load control~\cite{palensky2011demand} and our goal is to track a setpoint. Reference~\cite{kimonline} also provides a bandit extension to \texttt{COMID}. In our prior work~\cite{oco_irep2017}, we modified their bandit formulation for setpoint tracking. In this paper we also give two novel limited feedback extensions to \texttt{COMID}. Our contributions are:
\begin{itemize}
    \item We formulate an OCO-based setpoint tracking model for flexible loads (Section~\ref{sec:oco_formulation});
    \item We introduce a mean regularizer to minimize the impact of DR on loads;
    \item We propose a bandit-\texttt{COMID} algorithm for setpoint tracking (Section~\ref{sec:bandit});
    \item We introduce a partial bandit feedback extension in which only a subset of loads are observed, and provide a sublinear regret bound. (Section~\ref{ssec:pbf});
    \item We introduce a Bernoulli feedback extension in which the aggregator receives full feedback in some rounds and bandit feedback in others, and prove that the algorithm achieves sublinear regret (Section~\ref{ssec:berf});
    \item We numerically demonstrate the performance of our algorithms for setpoint tracking with TCLs and EVs. (Section~\ref{sec:num}).
\end{itemize}

\section{Background}
In this section, we introduce OCO and the \texttt{COMID} algorithm.

\subsection{Online convex optimization}
In OCO, a player chooses a decision from a convex set and suffers a loss, which is unknown at the time of the decision~\cite{hazan2016introduction}. At the end of each round, the player observes that round's loss. We index the rounds by $t$, denote the time horizon $T$ and the decision variable $\xx \in \mathcal{K}$. The performance of an online learning algorithm is characterized by its regret, defined as
\begin{equation*}
R_T = \sup_{\{F_1, F_2, \ldots, F_T \} \subset \mathcal{L}} \left\{ \sum_{t=1}^T F_t(\xx_t) - \min_{\xx \in \mathcal{K}} \sum_{t=1}^T F_t(\xx) \right\},
\end{equation*}
where $F_t$ is the loss function at time $t$ and $\mathcal{L}$ is the set of loss functions. $\mathcal{K} \subseteq \Re^N$ is the convex and compact decision set. The regret compares the cumulative loss suffered in each round by the player to the cumulative loss of the best fixed decision in hindsight. An algorithm that achieves a regret that is sublinear in the number of rounds eventually performs at least as well as the best fixed decision in hindsight~\cite{kimonline}.

\subsection{Composite objective mirror descent}
Online Gradient Descent (\texttt{OGD}) was first proposed in~\cite{zinkevich2003online} and then generalized to Online Mirror Descent (\texttt{OMD})~\cite{beck2003mirror}. In this work, we use the Composite Objective MIrror Descent (\texttt{COMID})~\cite{duchi2010composite}. This generalization of \texttt{OMD} handles loss functions of the form $F_t = f_t + r$, where $r$ is a round-independent regularizer. 

Define $\RR : \KK \rightarrow \Re$ as a regularization function (cf. references~\cite{hazan2016introduction,shalev2011online} for more detail about regularization functions). We only consider $\alpha$-strongly convex regularization functions $\RR$. The Bregman divergence with respect to $\RR$ is defined as
\begin{equation*}
\BR(\xx,\mathbf{z}) = \RR(\xx) + \RR(\mathbf{z}) - \nabla \RR (\mathbf{z})\Tr (\xx-\mathbf{z}).
\end{equation*}
Also define the dual norm of $\mathbf{z}$ as
\begin{equation*}
\norm{\mathbf{z}}_\ast = \sup\left\{\left. \mathbf{z}\Tr \xx  \ \right| \norm{\xx} \leq 1  \right\}.
\end{equation*}
Note that the dual of the $\ell_2$-norm is the $\ell_2$-norm.

The \texttt{COMID} update is given by
\begin{equation}
\xx_{t+1} = \argmin_{\xx \in \mathcal{K}} \eta \nabla_{\xx} f_t(\xx_t) \Tr \xx + \BR(\xx, \xx_t) + \eta r(\xx), \label{eq:update}
\end{equation}
where $\eta$ is a numerical parameter. We now present a specialized version of the original \texttt{COMID} regret bound of~\cite{duchi2010composite}. In this version, $\eta$ is tuned to avoid too small step sizes when applied to setpoint tracking.

\begin{lemma}[Regret bound for \texttt{COMID}]
Let $f_t$ be a $L$-Lipschitz function and $\norm{\nabla f_t(\mub_t)}_\ast \leq G_\ast$ for all $t$, $r(\mub_1) = 0$ and define the tuning parameter $\chi \geq 1$. Then, using
\begin{equation}
\eta = \chi \sqrt{\frac{2 \alpha \BR(\mub^\ast, \mub_1)}{G_\ast^2 T}}, \label{eq:eta_duchi}
\end{equation}
the regret of \texttt{COMID} is upperbounded by,
\begin{equation}
R_T(\texttt{COMID}) \leq \sqrt{\frac{2T \BR(\mub^\ast, \mub_1)G_\ast^2 \chi^2}{\alpha}}
\end{equation}
\label{lem:tun-comid}
\end{lemma}
The proof of Lemma~\ref{lem:tun-comid} can be found in the conference version of this paper~\cite[Lemma 1]{oco_irep2017} and relies on~\cite[Corollary 4]{duchi2010composite}.

In this work, we set $\RR(\cdot) = \frac{1}{2}\| \cdot \|^2_2$. The Bregman divergence $\BR$ then simplifies to $\BR(\xx,\mathbf{z}) = \frac{1}{2}\norm{\xx - \mathbf{z}}^2_2$, a $1$-strongly convex function. Under this regularization function, the \texttt{OMD} algorithm simplifies to \texttt{OGD}. Accordingly we refer to \texttt{COMID} with $\RR(\cdot) = \frac{1}{2}\| \cdot \|^2_2$ as the Composite Objective Gradient Descent (\texttt{COGD}). For the rest of this work, we will be using the \texttt{COGD} algorithm. Finally, we let $D$ be the diameter of the compact set $\KK$, 
\[
D = \mathrm{diam} \ \mathcal{K} = \sup\left\{ \left. \norm{\xx - \mathbf{z}}_2 \, \right| \, \xx, \mathbf{z} \in \mathcal{K} \right\},
\]
and we note that $\BR(\xx,\mathbf{z}) \leq \frac{1}{2} D^2$ for $\xx, \mathbf{z} \in \mathcal{K}$.

\section{Setpoint tracking with online convex optimization}
\label{sec:oco_formulation}
The objective of the load aggregator is to send adjustment signals to loads so that the total adjustment tracks a power setpoint, $s_t \in \Re$. We consider $N$ flexible loads. We let $\mub_t \in \left[-1, 1\right]^N$ denote the vector of adjustment signals sent to the loads at time $t$. $\mub_t$ represents instructions to scale down ($ \mu < 0$) or scale up ($ \mu > 0$) the load power consumption to match the setpoint. For example, in the case of Thermostatically Controlled Loads (TCLs), scaling down its power consumption represents lowering the air conditioning intensity or simply shutting it down, and scaling up represents increasing intensity to stock thermal energy. This signal is scaled by the response of the loads, $\cc_t$, and the sum over all the loads gives the aggregate power adjustment. Depending on the feedback structure, the aggregate adjustment and part of $\cc_t$ is revealed immediately after the decision in the current round. The uncertainty in $\cc_t$ represents, for instance, erratic consumer behavior or the loads' dependence on weather.

We define the setpoint tracking loss as the square of the tracking error:
\begin{equation}
\ell(\mub_t) = (s_t - \cc \Tr_t \mub_t)^2. \label{eq:ell_loss}
\end{equation}
Note that $s_t$ can be known or unknown when $\mub_t$ is chosen.

We also use a mean-regularizer to penalize deviations of the mean adjustment over time from zero.
The purpose of this is to minimize the impact of demand response on the loads, e.g., to reduce  rebound effects~\cite{palensky2011demand}. Let $\langle \cdot \rangle_t$ denote the mean of its argument over rounds $1$ to $t$ and let the mean of a vector be the vector of the means of its elements. We define the mean-regularizer as
\begin{equation*}
\norm{\langle \mub \rangle_t}^2_2 = \norm{\frac{1}{t} \sum_{s=1}^t \mub_s}^2_2 =
\norm{\frac{(t-1) \langle \mub \rangle_{t-1} + \mub_t}{t}}^2_2.
\end{equation*}

We use an $\ell_1$-norm sparsity regularizer to minimize the number of loads dispatched in each round and to avoid dispatching loads with small contributions.

The total objective function is therefore given by
\begin{equation}
F_t(\mub_t) =\left(s_t - \cc\Tr_t \mub_t \right)^2 + \rho \norm{\langle \mub \rangle_t}^2_2 + \lambda \norm{\mub_t}_1, \label{eq:total_loss} 
\end{equation}
where $\rho$ and $\lambda$ are numerical parameters.

We now apply COGD to setpoint tracking. We define the loss function and regularizer to be
\begin{align}
f_t(\mub_t) &= \left(s_t - \cc\Tr_t \mub_t \right)^2 + \rho \norm{\langle \mub \rangle_t}^2_2 \label{eq:f_comid},\\
r(\mub_t) &= \lambda \norm{\mub_t}_1. \label{eq:r_comid}
\end{align}
Note that the mean-regularizer is included in~\eqref{eq:f_comid} and not in~\eqref{eq:r_comid} because \texttt{COGD}, like the more general \texttt{COMID}, only handles round-independent regularizers.
The adjustment signal update is obtained applying~\eqref{eq:total_loss} to~\eqref{eq:update}, which yields
\begin{align}
\mub_{t+1} &= \argmin_{\mub \in [-1,1]^N} \left\{ \vphantom{\left[-2 \cc_t (s_t - \cc_t \Tr \mub_t)  + \frac{2\rho}{t} \frac{(t-1) \langle \mub \rangle_{t-1} + \mub_t}{t}\right]\Tr} \frac{1}{2}\norm{\mub_t - \mub }_2^2 + \eta \lambda \norm{\mub}_1 \right. \label{eq:update_comid}\\
+& \left. \eta \left[-2 \cc_t (s_t - \cc_t \Tr \mub_t)  + \frac{2\rho}{t} \frac{(t-1) \langle \mub \rangle_{t-1} + \mub_t}{t}\right] \Tr \mub \right\}. \nonumber
\end{align}
Our \texttt{COGD} algorithm for setpoint tracking is given in a preliminary version of this work~\cite[Figure 1]{oco_irep2017}. We conclude this section by presenting Proposition~\ref{prop:comid}. This result follows from Lemma~\ref{lem:tun-comid} and establishes that our \texttt{COGD} algorithm for setpoint tracking achieves sublinear regret.

\begin{proposition}
\label{prop:comid}
Let $\mub_1 = \mathbf{0}$, $\chi \geq 1$, $f_t(\mub_t) \leq B$ for all $t$ and set
\begin{equation}
\eta = \chi \sqrt{\frac{4N}{G^2T}} \label{eq:eta_alg}.
\end{equation}
Then, the \texttt{COGD} for setpoint tracking algorithm has a regret bound given by
\begin{equation}
R_T\leq 4\chi\sqrt{T K B}, \label{eq:regret_alg1}
\end{equation}
where $\displaystyle K = \max_{t=1,2, \ldots, T}\{ \rho^2, \norm{\cc_t}^2_2\}$.
\end{proposition}
The reader is referred to~\cite{oco_irep2017} for the proof of Proposition~\ref{prop:comid}.

We subsequently refer to the assumptions of this section as the full information setting.

\section{Limited feedback}
\label{sec:lim}
In this section, we present three limited-feedback extensions to the OCO model introduced in the last section. We give bounded-regret algorithms for each extension. 

\subsection{Bandit feedback}
\label{sec:bandit}

In the bandit feedback setting, the aggregator only observes aggregate effect of its decision in each round. Specifically, it sees $\cc_t \Tr \mub_t$ and $s_t$. Hence, no gradient can be computed for the \texttt{COGD} update~\eqref{eq:update_comid}. Using only the total loss, a point-wise gradient estimator $\gt$ can be computed at each round. We refer the reader to~\cite{oco_irep2017} for the detail derivation of the point-wise gradient, which is based on the approach of~\cite{flaxman2005online,hazan2016introduction}.

We define the point-wise gradient estimator for setpoint tracking as
\begin{equation}
\gt = \frac{N}{\delta} f_t(\mub_t + \delta \vb_t) \vb_t,
\label{eq:est}
\end{equation}
where $\delta > 0$ is a numerical parameter, $\vb_t$ is sampled uniformly from $\Sone^N=\left\{\left. \vb \in \Re^N \right| \norm{\vb}_2 = 1 \right\}$ and $f_t$ is given in~\eqref{eq:f_comid}. The expectation of the point-wise gradient~\eqref{eq:est} over $\vb_t$ is equal to the gradient of $f_t$.

We must also modify the decision set so that the estimated gradient does not lead to an infeasible step. Define the convex set $\mathcal{K}^\delta \subset \mathcal{K}$ as
\begin{align}
\mathcal{K}^{\delta} &\equiv \left\{ \mub \left| \frac{\mub}{1-\delta} \in \mathcal{K} \right. \right\} = \left[\delta-1, 1-\delta \right]^N. \label{eq:delta_set}
\end{align}

The update for the bandit-\texttt{COGD} (\texttt{BCOGD}) algorithm is then given by
\begin{equation}
\mub_{t+1} = \argmin_{\mub \in \mathcal{K}^{\delta}}  \eta \gt \Tr \mub + \frac{1}{2}\norm{\mub_t - \mub }_2^2 + \eta r(\mub). \label{eq:update_bandit_cogp}
\end{equation}
The full algorithm is given in our preliminary work~\cite[Figure 2]{oco_irep2017}.
\begin{theorem}[Regret bound for \texttt{BCOGD}]
\label{thm:bcogd}
Let $F_t(\mub_t)$ be $L$-Lipschitz and $B$-bounded for all $t$ and let $r(\mub_1)=0$. Then, using the point-wise gradient estimator $\gt$ and setting
\begin{equation}
\eta = \frac{D \chi}{B N T^{\frac{3}{4}}},  \qquad\qquad \delta = \frac{1}{T^{\frac{1}{4}}} 
\label{eq:eta_bcogp}
\end{equation}
where $D = \mathrm{diam} \ \mathcal{K}$ and $\chi \geq 1$, the \texttt{BCOGD} regret is upper bounded by
\begin{equation}
\E[R_T(\texttt{BOGD})] \leq \left(DBN\chi + 2DL + 2L \right) T^{\frac{3}{4}}. \label{eq:bcogd_bound}
\end{equation}
\end{theorem}
The proof of Theorem~\ref{thm:bcogd} is given in~\cite{oco_irep2017}. The $O\left(T^{3/4}\right)$ bound is similar to the regret of other bandit OCO algorithms~\cite{flaxman2005online,hazan2016introduction,kimonline}. Observe that there is an increase in the round-dependence of the bound from $1/2$ to $3/4$ due to less information being available to the algorithm.

This result is a generalization of the bandit-\texttt{COGD} proposed in~\cite{kimonline}, which is based on the bandit \texttt{OGD} of~\cite{flaxman2005online}. Our generalization solves the issue of~\cite{kimonline} requiring a very large time horizon and small $\eta$, the update step size, in the case of setpoint tracking. A more detailed discussion of these issues is given in~\cite{oco_irep2017}.
 
When applied to setpoint tracking, the regret of \texttt{BCOGD} reduces to
\begin{equation*}
\E\left[R_T(\texttt{BCOGD})\right] \leq \left(2N^{\frac{3}{2}}B\chi + 4 \sqrt{N}L\right)T^{\frac{3}{4}}.
\end{equation*}

\subsection{Partial Bandit Feedback} 
\label{ssec:pbf}
We now present the Partial Bandit \texttt{COGD} (\texttt{PBCOGD}) algorithm for when the aggregator only receives full feedback from a subset of the loads and has access to the aggregate effect of its decision, $\cc \Tr \mub$. This could represent a total power measurement, e.g., at a substation. This extension applies, for example, when individual monitoring has only been implemented for a subset of the loads, when some loads have opted out due to privacy issues, or when the aggregator is subject to bandwidth constraints. We drop the mean regularizer in \texttt{PBCOGD} because we have observed empirically that it has little effect in this setting.

Let $n \in \{1,2, \ldots, N-1 \}$ be the number of loads for which full information is available in each round. 

We denote $\mub_t^F \in \mathcal{K}_n$ and $\mub_t^B \in\mathcal{K}_{N-n}$ as the decisions variables and sets for the full information and the bandit-like load subsets, where $\mathcal{K}_n \subseteq \Re^n$ and $\mathcal{K}_{N-n} \subseteq \Re^{N-n}$ are compact and convex sets. We can write $\mub_t$ as $\mub_t = \tilde{\mub}_t^{B} + \tilde{\mub}_t^{F}$ where $\tilde{\mub}_t^{B} = \left(\mub_t^B, \mathbf{0} \right)\Tr$ and $\tilde{\mub}_t^{F} = \left(\mathbf{0},\mub_t^F \right)\Tr$. 

Define $\beta_t=\cc_t \Tr \mub_t - \cc_t^{F \, \mathrm{T}} \mub_t^F$, where $\cc_t^F$ is the subset of responses corresponding to the full information loads. The update for the \texttt{PBCOGD} is,
\begin{align}
\mub_{t+1} =& \argmin_{\substack{\mub^F \in \mathcal{K}_n \\ \mub^B \in \mathcal{K}_{N-n}^{\delta}}}  \Bigg\{\frac{1}{2}\norm{\mub_t - \mub }_2^2 + \eta_1 r(\mub^B) + \eta_2 r(\mub^F) \nonumber\\
&+ \eta_1 \gt^{p \, \mathrm{T}} \mub_t^B + \eta_2 \nabla_{\mub^F} \tilde{f}_t^F(\beta_t, \mub_t^F) \Tr \mub_t^F \Bigg\} \label{eq:update_partial}
\end{align}
where $\tilde{f}_t^F(\beta_t, \mub_t^F)$ is a function of $\beta_t$ and the full information decision variables and is equal to $f_t(\mub_t)$. The definition of $\gt^p$, the point-wise gradient for \texttt{PBCOGD}, is given with the full algorithm in Figure~\ref{alg:partial}. We then have the following regret bound.

\begin{theorem}[Regret bound for \texttt{PBCOGD}]
\label{thm:partial}
Assume that $r(\mub_t)$, the regularizer, satisfies $r(\mub_t) = r(\mub_t^B) + r(\mub_t^F)$ and $\beta_t$ is available for all $t$. Let $\eta_1$ and $\delta$ be defined as in~\eqref{eq:eta_bcogp} and set $\eta_2$ as in~\eqref{eq:eta_duchi}. 
Then, under the assumptions of Lemma~\ref{lem:tun-comid} and Theorem~\ref{thm:bcogd}, \texttt{PBCOGP} achieves $O\left( T^{\frac{3}{4}} \right)$ regret bound.
\end{theorem}
The proof of Theorem~\ref{thm:partial} is given in Appendix~\ref{app:partial}.

We now apply the \texttt{PBCOGD} to setpoint tracking. We have $r(\mub) = \lambda \norm{\mub}_1$ which respects the regularizer condition.

We denote $\cc_t^B$ and $\cc_t^F$ as the subset of response corresponding to the bandit and full information loads, respectively. We have
\begin{equation*}
\tilde{f}_t^F(\mub_t^F) = \left(s_t - \beta_t - \cc\Tr_t \begin{pmatrix} \mathbf{0} \\ \mub_t^F \end{pmatrix} \right)^2.
\end{equation*}
The \texttt{PBCOGD} for setpoint tracking with partial feedback algorithm is presented in Figure~\ref{alg:partial}. By Theorem~\ref{thm:partial}, the regret of \texttt{PBCOGD} for setpoint tracking is sublinear.

\begin{figure}
\fbox{
\begin{minipage}[b]{0.94\columnwidth}
\begin{algorithmic}[1]
\STATE \textbf{Parameters:} Given $n$, $T$, $\rho$ and $\chi$.
\STATE \textbf{Initialization:} Set $\mub_1=\mathbf{0}$ and set $\eta_1$ and $\delta$ according to~\eqref{eq:eta_bcogp} and $\eta_2$ according to \eqref{eq:eta_duchi}.
\medskip

\FOR{$t = 1,2, \ldots, T$}
\STATE Sample $\vb^p_t \sim \Sone^{N-n}$.
\STATE Dispatch adjustment according to $$\mub_t = \left[ \mub_t^B + \delta \vb^p_t , \mub_t^F \right]\Tr.$$
\STATE Suffer loss $\ell_t(\mub_t)$ and compute 
\[
\beta_t = \cc_t \Tr \mub_t - \cc_t^{F \, \mathrm{T}} \mub_t^F .
\]
\STATE Compute the gradients,
\begin{equation*}
\gt^p = \frac{(N-n)}{\delta} f_t\left(\mub_t + \delta \begin{pmatrix} \vb^p_t \\ \mathbf{0} \end{pmatrix}\right) \ \mathbf{v}^p_t,
\end{equation*}
\begin{align*}
\nabla_{\mub^F} \tilde{f}_t^F(\beta_t, \mub_t^F) &= -2 \cc_t^F (s_t - \beta_t - \cc_t^{F \, \mathrm{T}} \mub_t^F).  \\
\end{align*}
\STATE Update adjustment dispatch,
\begin{align*}
\mub_{t+1} &= \begin{bmatrix}
\mub^B \\ \mub^F
\end{bmatrix} \\
&= \argmin_{\substack{\mub^F \in \mathcal{K}_n \\ \mub^B \in \mathcal{K}_{N-n}^{\delta}}}  \Big\{ \eta_1 \gt^{p \, \mathrm{T}} \mub_t^B  \\
&\quad + \eta_2 \nabla_{\mub^F} \tilde{f}_t^F(\mub_t^F) \Tr \mub_t^F + \eta_1 \lambda \norm{\mub^B}_1 \\
& \quad +\eta_2 \lambda \norm{\mub^F}_1 + \frac{1}{2}\norm{\begin{pmatrix} \mub^B_t \\ \mub^F_t\end{pmatrix} - \begin{pmatrix} \mub^B \\ \mub^F\end{pmatrix} }_2^2 \Big\}
\end{align*}
\ENDFOR
\end{algorithmic}
\end{minipage}
}
\caption{\texttt{PBCOGB} for setpoint tracking with partial feedback algorithm}
\label{alg:partial}
\end{figure}

\subsection{Bernoulli feedback setting}
\label{ssec:berf}
We now consider an extension in which the decision maker receives full feedback in some rounds and bandit feedback in others, depending on the outcome of a Bernoulli random variable. This decreases the communication requirement. During the initialization of the algorithm, the player randomly determines the feedback types of all rounds. We present here the Bernoulli-\texttt{COGD} (\texttt{BerCOGD}) algorithm to deal with this new type of feedback. 

If the probability of bandit feedback in a given round is below a threshold, we can obtain a regret upper-bound of $O\left(T^{1/2}\right)$ for the \texttt{BerCOGD} algorithm, similarly to the \texttt{COGD} algorithm for full feedback. We denote the feedback type in each round with the random variable $I_t \sim \mathrm{Bernoulli}(p)$ for $t=1,2, \ldots, T$, with $I_t = 1$ referring to a bandit feedback and $I_t=0$ to full feedback.

\begin{figure}
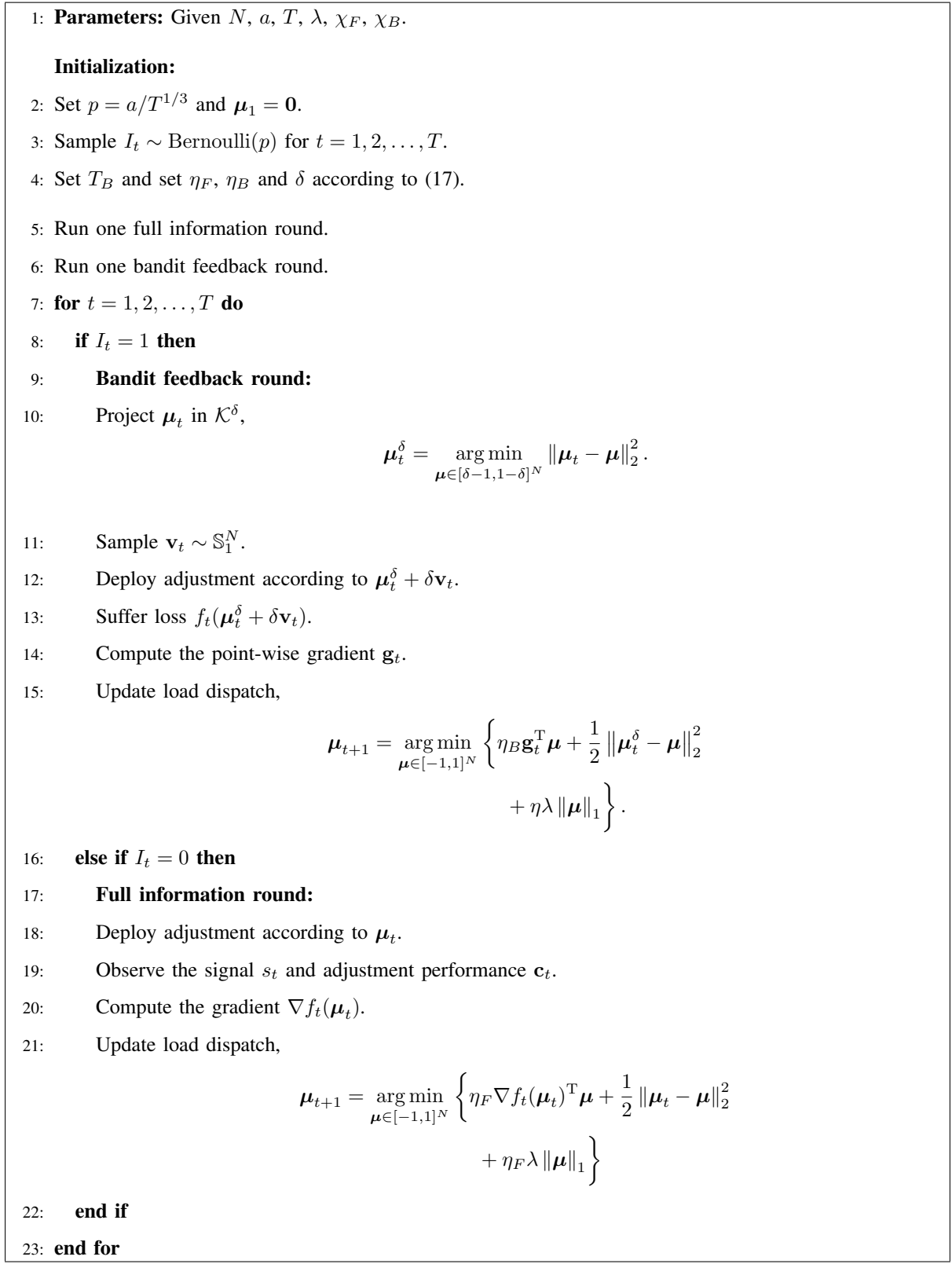

\fbox{
\begin{minipage}[b]{0.94\columnwidth}
\begin{algorithmic}[1]
\STATE{\textbf{Parameters:} Given $N$, $a$, $T$, $\lambda$, $\chi_F$, $\chi_B$.\\}
\medskip

\textbf{Initialization:} 
\STATE{Set $p = a/T^{1/3}$ and $\mub_1 = \mathbf{0}$.}
\STATE{Sample $I_t \sim \mathrm{Bernoulli}(p)$ for $t = 1, 2, \ldots, T$.}
\STATE{ Set $T_B$ and set $\eta_F$, $\eta_B$ and $\delta$ according to~\eqref{eq:etaF_ber-cogp}.\\}
\medskip

\STATE{Run one full information round.}
\STATE{Run one bandit feedback round.}

\FOR{$t = 1,2, \ldots, T$}
\IF {$I_t = 1$}
\STATE \textbf{Bandit feedback round:}
\STATE{Project $\mub_t$ in $\mathcal{K}^{\delta}$,
\[
\mub_t^{\delta} = \argmin_{\mub \in [\delta -1 , 1- \delta]^N} \norm{\mub_t - \mub}^2_2.
\]}
\STATE{Sample $\vb_t \sim \Sone^N$.}
\STATE {Deploy adjustment according to $\mub_t^{\delta} + \delta \vb_t$.}
\STATE {Suffer loss $f_t(\mub_t^{\delta} + \delta \vb_t)$.}
\STATE {Compute the point-wise gradient $\gt$.
}
\STATE {Update load dispatch,
\begin{align*}
\mub_{t+1} = \argmin_{\mub \in [-1,1]^N}  &\left\{\eta_B \gt \Tr \mub + \frac{1}{2}\norm{\mub_t^{\delta} - \mub }_2^2 \right.
\\ &\left. \vphantom{\frac{1}{2}} \quad+ \eta \lambda \norm{\mub}_1 \right\}.
\end{align*}}
\ELSIF {$I_t = 0$}
\STATE \textbf{Full information round:}
\STATE {Deploy adjustment according to $\mub_t$.}
\STATE {Observe the signal $s_t$ and adjustment performance $\cc_t$.}
\STATE {Compute the gradient $\nabla f_t(\mub_t)$.
}
\STATE {Update load dispatch,
\begin{align*}
\mub_{t+1} = \argmin_{\mub \in [-1,1]^N}  &\left\{\eta_F \nabla f_t(\mub_t) \Tr \mub + \frac{1}{2}\norm{\mub_t - \mub }_2^2 \right.\\
& \left.\vphantom{\frac{1}{2}} \quad+\eta_F \lambda \norm{\mub}_1 \right\}
\end{align*}
}
\ENDIF{}
\ENDFOR
\end{algorithmic}
\end{minipage}
}
\caption{\texttt{BerCOGD} for setpoint tracking with Bernoulli feedback algorithm}
\label{alg:bernoulli}
\end{figure}

The \texttt{BerCOGD} is as follows: if the $I_t =0$ at round $t$, then the \texttt{COGD} update is used and if $I_t =1$, then the \texttt{BCOGD} update is used. The only difference is that in the bandit feedback update the projection of $\mub_t$ onto $\KK^\delta$ is done in two steps.
When the decision $\mub_{t+1}$ using the bandit update is computed, it is projected onto $\KK$ rather than $\KK^\delta$ as shown on line 15 of the algorithm in Figure~\ref{alg:bernoulli}. The projection onto $\KK^\delta$ is made at the beginning of a bandit round (cf. line 10 of the algorithm). This is necessary to ensure that a bandit round can follow a full information round.

\begin{theorem}[Regret bound for \texttt{BerCOGD}]
Let the feedback type in round $t$ be determined by the random variable $I_t\sim \mathrm{Bernoulli}(p)$ where $p\leq \frac{a}{T^{1/3}} \leq 1$. Then, under the assumptions of Lemma~\ref{lem:tun-comid} and Theorem~\ref{thm:bcogd} and setting
\begin{align}
\eta_F =  \frac{D \chi_F}{G \left(T-T_B + 1\right)^\frac{1}{2}}&, \  \eta_B = \frac{D \chi_B}{B N (T_B+1)^{\frac{3}{4}}}, \nonumber\\
\delta =& \frac{1}{(T_B+1)^{\frac{1}{4}}}, \label{eq:etaF_ber-cogp}
\end{align}
the \texttt{BerCOGD} algorithm's expected regret is bounded by,
\begin{equation}
\E[R_T({\texttt{BerCOGD}})] \leq \left(A_1 + A_2 a^{\frac{3}{4}}\right) T^{\frac{1}{2}} + A_1 + A_2, \label{eq:bern_regret_gen}
\end{equation}
where $A_1 = DG\chi_F$, $A_2 = DBN\chi_B + 2DL + 2L$. 
\label{thm:bernoulli}
\end{theorem}

The proof of Theorem~\ref{thm:bernoulli} is given in Appendix~\ref{app:bern}. We note that as $T$ gets larger, we can accommodate less and less bandit rounds. However, the constant $a$ of the probability $p$ can be tuned to increase the ratio of bandit rounds. Another approach is to choose a smaller $T$ and periodically reinitialize the algorithm.

We now apply \texttt{BerCOGD} to the setpoint tracking problem. The algorithm is shown in Figure~\ref{alg:bernoulli}. The gradient $\mathbf{w}_t$ and point-wise gradient $\gt$ are define as in Section~\ref{sec:oco_formulation} and Section~\ref{sec:bandit} but without the mean-regularizer term. Applying Theorem~\ref{thm:bernoulli} to the \texttt{BerCOGD} for setpoint tracking algorithm, the regret bound~\eqref{eq:bern_regret_gen} becomes
\begin{align}
\E[R_T] &\leq \left(\sqrt{16KB\chi_F^2} + DBN\chi_Ba^{\frac{3}{4}} +2DLa^{\frac{3}{4}} \right) T^{\frac{1}{2}}\nonumber\\
&\qquad + \sqrt{16KB\chi_F^2} + \left(DBN\chi_B +2DL\right) \label{eq:bern_regret}
\end{align}
Thus,~\eqref{eq:bern_regret} shows that the regret upper bound is $O(\sqrt{T})$ as in the full information setting even if on average $Tp$ rounds only have access to bandit feedback.

\section{Numerical results}
\label{sec:num}
We now apply our OCO algorithms to setpoint tracking with TCLs.

\subsection{Modeling}
\label{sec:tcls}
We use the model developed by~\cite{mortensen1988stochastic,ucak1998effects} and used in DR by~\cite{callaway2009tapping,mathieu2015arbitraging} to model the temperature dynamics of a TCL. The temperature of a TCL, $\theta_t$, for a time step of length $h$ evolves as
\begin{equation}
\theta_{t+1} = b \theta_{t} + (1-b) \left( \theta_a - m R P_r \right), \label{eq:temp_dyn}
\end{equation}
where $\theta_a$ is the ambient temperature, $b=\ee^{-h/RC}$, $m \in \{0,1\}$ is the cooling unit control and $R$, $C$, $P_R$ are thermal parameters. Note that this model also applies to heating. 

To make the model amenable to OCO, we relax the binary control to $m\in[0,1]$, a common assumption~\cite{nguyen2014optimal,nguyen2014joint}. We assume that the unit tries to stay at its desired temperature $\theta_d$ by setting its cooling unit control to $m = \overline{m}$. Hence, by fixing $\theta_a$ for all $t$, we define $\overline{m}$ as
\[
\overline{m} = \frac{\theta_a - \theta_d}{P_R R}.
\]
Any deviation from $\overline{m}$ will represent an increase or decrease in power consumption and can hence be used as an adjustment for setpoint tracking. The power consumption of the TCLs can be set to any value in the interval $[-\overline{m}, 1-\overline{m}]$. We constrain the adjustment interval to be symmetric so the maximum increase or decrease in power consumption is equal in absolute value.
For all loads $i$, we define the average adjustment response as
\[
c_0(i) =p(i) \min \left\{\overline{m}(i),1-\overline{m}(i) \right\} 
\]
where $p(i)=P_R(i)/COP(i)$. The cooling unit control is
\[
m(i) = \mu(i) \min \left\{\overline{m}(i),1-\overline{m}(i) \right\} .
\]
In the definition of $c_0(i)$, $P_R$ is the rated power and $COP$ is the coefficient of performance of the cooling unit. 

We rewrite the setpoint tracking loss function~\eqref{eq:ell_loss} as
\begin{equation*}
\ell(\mub_t) = (s_t - \mathbf{p}_t \Tr \mathbf{\overline{m}} - \mathbf{c}_t \Tr \mub_t)^2.
\end{equation*}
The mean and sparsity regularizers are as in~\eqref{eq:total_loss}, and in this context penalize discomfort and using too many TCLs, respectively.

In our numerical simulations, we consider $N=100$ loads and let $n=10$ in the partial bandit feedback simulations. We set the ambient temperature to 30$^\circ$C and the time step to $h=5$ minutes. We set $s_t = 15 \sin (0.1 t) + 155$. 
We define the loads' response at round $t$ to be $c_t(i) = c_0(i) + w_t$, where $w_t\sim \mathrm{N}_{[-1,1]}(0,0.5)$, a truncated Gaussian variable, is used to model the uncertainty of the loads' response for all load $i$ and $t=1,2, \ldots T$.

The TCLs parameters are sampled uniformly from thermal parameters in~\cite{mathieu2012state} except for the desired temperature $\theta_d$, which is sample uniformly between 20$^{\circ}\text{C}$ and 25$^{\circ}\text{C}$ for all loads for the purpose of increasing the flexibility of the loads.

\subsection{Numerical results}
We now present the simulation results. The optimization problem in each update is solved numerically using \texttt{CVX}~\cite{cvx,gb08} and Gurobi~\cite{gurobi}. For each case except full information, the average over 100 trials is shown. The parameter $\chi$ is set to $200$ when \texttt{COGD} updates are used, to $5.5 \times 10^4$ when \texttt{BCOGD} updates are used, and to $150$ and $3 \times 10^4$ when respectively \texttt{COGD} and \texttt{BCOGD}  updates are used in the \texttt{BerCOGD} algorithm. The probability parameter $a$ is fixed to $7.6$ for the \texttt{BerCOGD} algorithm. The regularization parameters $\rho$ and $\lambda$ are given in Tables~\ref{tab:set_comp} and~\ref{tab:reg_comp} for each algorithm.

\begin{figure}[!t]
    \centering
    \includegraphics[width=0.75\columnwidth,height=0.75\textheight,keepaspectratio]{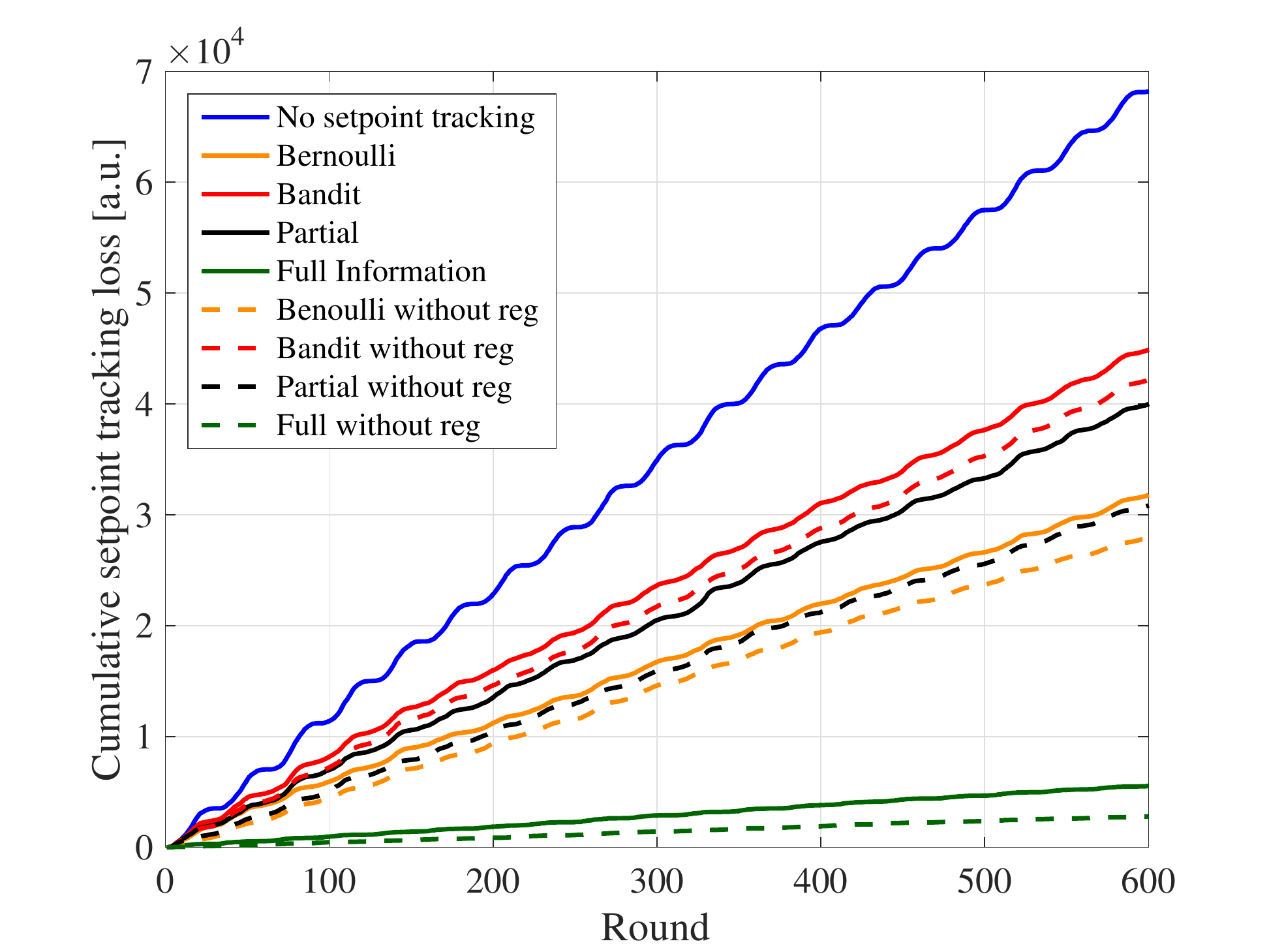}
    \caption{Cumulative setpoint tracking loss comparison through different feedback type}
    \label{fig:setpoint_loss}
\end{figure}

Figure~\ref{fig:setpoint_loss} presents the cumulative setpoint tracking loss for each algorithm with and without regularization. Table~\ref{tab:set_comp} shows the percentage improvement in the loss function compared to when no DR is performed. Table~\ref{tab:reg_comp} shows the impact of the regularizers on each algorithm. To quantify the effects of the mean and sparsity regularizers, we compare the Euclidean norm of the mean of $\mub_t$ up to the current round and the $1$-norm of $\mub_t$ with and without regularization. We average the result over all rounds. We note that the sparsity regularizer also promotes a zero mean and that the mean regularizer has a reduced impact on the total setpoint tracking loss while also promoting sparsity.

\begin{table}[!t]
\renewcommand{\arraystretch}{1.3}
\caption{Total setpoint tracking improvement comparison (averaged over 100 simulations)}
\label{tab:set_comp}
\centering
\begin{tabular}{c l l}
\hline
Feedback type & with regularization & without\\
\hline
Full information & $91.87\%$ ($\rho = 250$ \& $\lambda = 7.5$)& $95.89\%$\\
Bandit & $34.15\%$ ($\rho = 1.5$ \& $\lambda = 60$)& $38.12\%$\\
Partial Bandit ($n=10$) & $41.33\%$ ($\rho = 0$ \& $\lambda = 40$)& $54.74\%$\\
Bernoulli ($p=0.9$) & $53.39\%$ ($\rho = 2.5$ \& $\lambda = 65 $)& $58.96\%$\\
\hline
\end{tabular}
\end{table}

\begin{table}[!t]
\renewcommand{\arraystretch}{1.3}
\caption{Per round average regularizer improvement (averaged over 100 simulations)}
\label{tab:reg_comp}
\centering
\begin{tabular}{c l l l}
\hline
Feedback type & Mean &  Sparsity & Parameters\\
\hline
Full information & $77.90\%$ & $34.15\%$ & $\rho = 250$ \& $\lambda = 7.5$ \\
Bandit & $25.72\%$ & $5.29\%$ & $\rho = 1.5$ \& $\lambda = 60$\\
Partial Bandit (n=10) & \begin{centering} ---\end{centering} & $5.70\%$ & $\lambda = 40$\\
Bernoulli ($p=0.9$) & $52.57\%$ & $25.03\%$ & $\rho = 2.5$ \& $\lambda = 65 $ \\
\hline
\end{tabular}
\end{table}
Figure~\ref{fig:setpoint_loss} and Table~\ref{tab:set_comp} show how feedback can improve the setpoint tracking loss. The full information algorithm performs significantly better than the bandit and partial feedback algorithms. The Bernoulli feedback extension, however, offers the second best performance while receiving on average full feedback in $10\%$ of the rounds and bandit in the rest. 
Lastly, we compare the bandit and the partial bandit algorithms. Figure~\ref{fig:setpoint_loss} shows the improvement of the cumulative loss by the partial bandit when only $10\%$ of full information variables are available to the load aggregator.
 
\begin{figure}[!t]
\centering
\includegraphics[width=0.75\columnwidth,height=0.75\textheight,keepaspectratio]{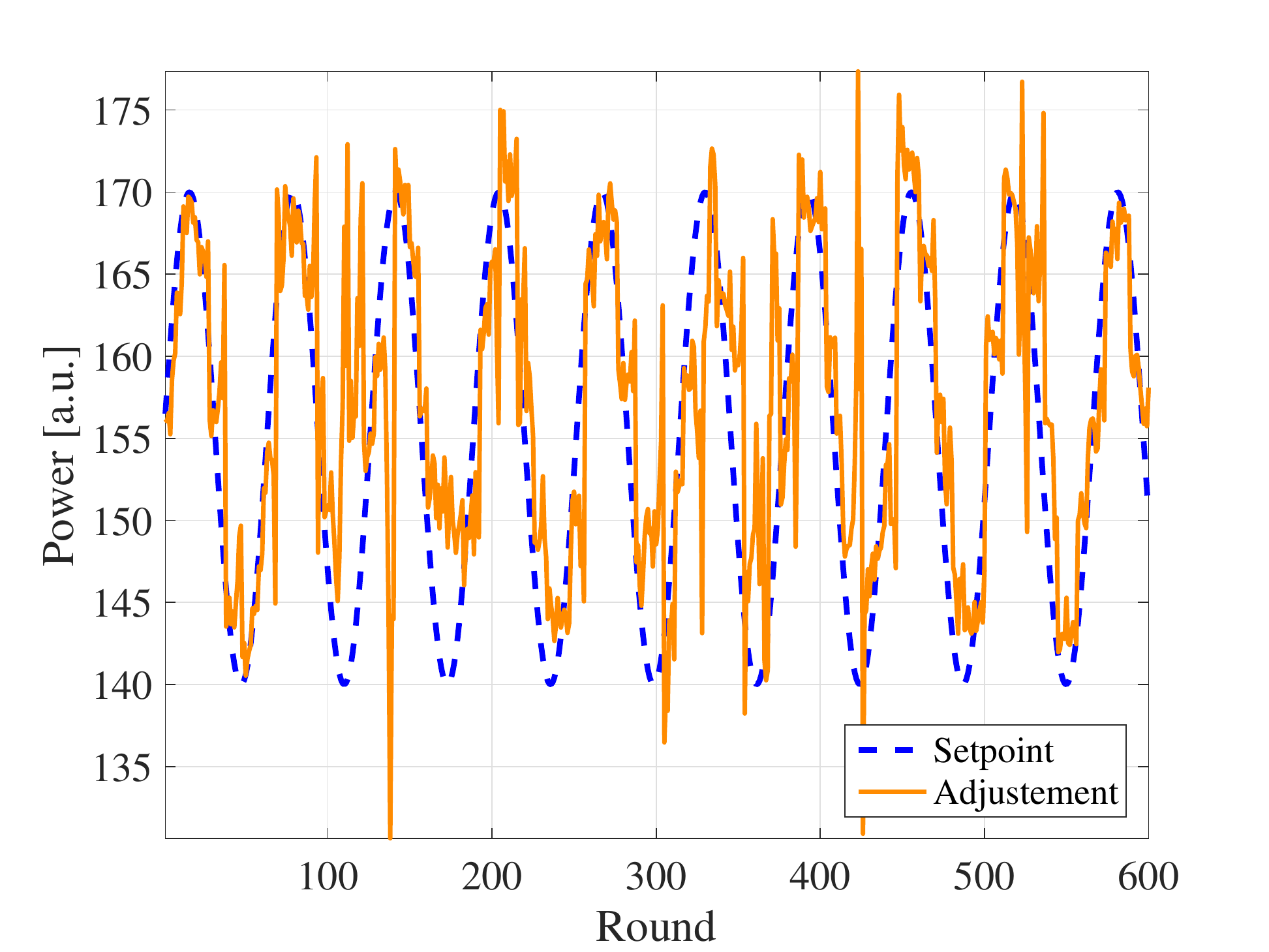}
\caption{Setpoint tracking under Bernoulli feedback}
\label{fig:setpoint_partial}
\end{figure}

An instance of the setpoint tracking curve is presented in Figure~\ref{fig:setpoint_partial} for the Bernoulli feedback algorithm. Setpoint tracking curves for the full information (with and without regularization) and bandit feedback can be found in~\cite{oco_irep2017}. 

\begin{figure}[!t]
\centering
\subfloat[Non-zero signal rounds for selected TCLs (tolerance $10^{-2}$)]{\includegraphics[width=0.75\columnwidth,height=0.75\textheight,keepaspectratio]{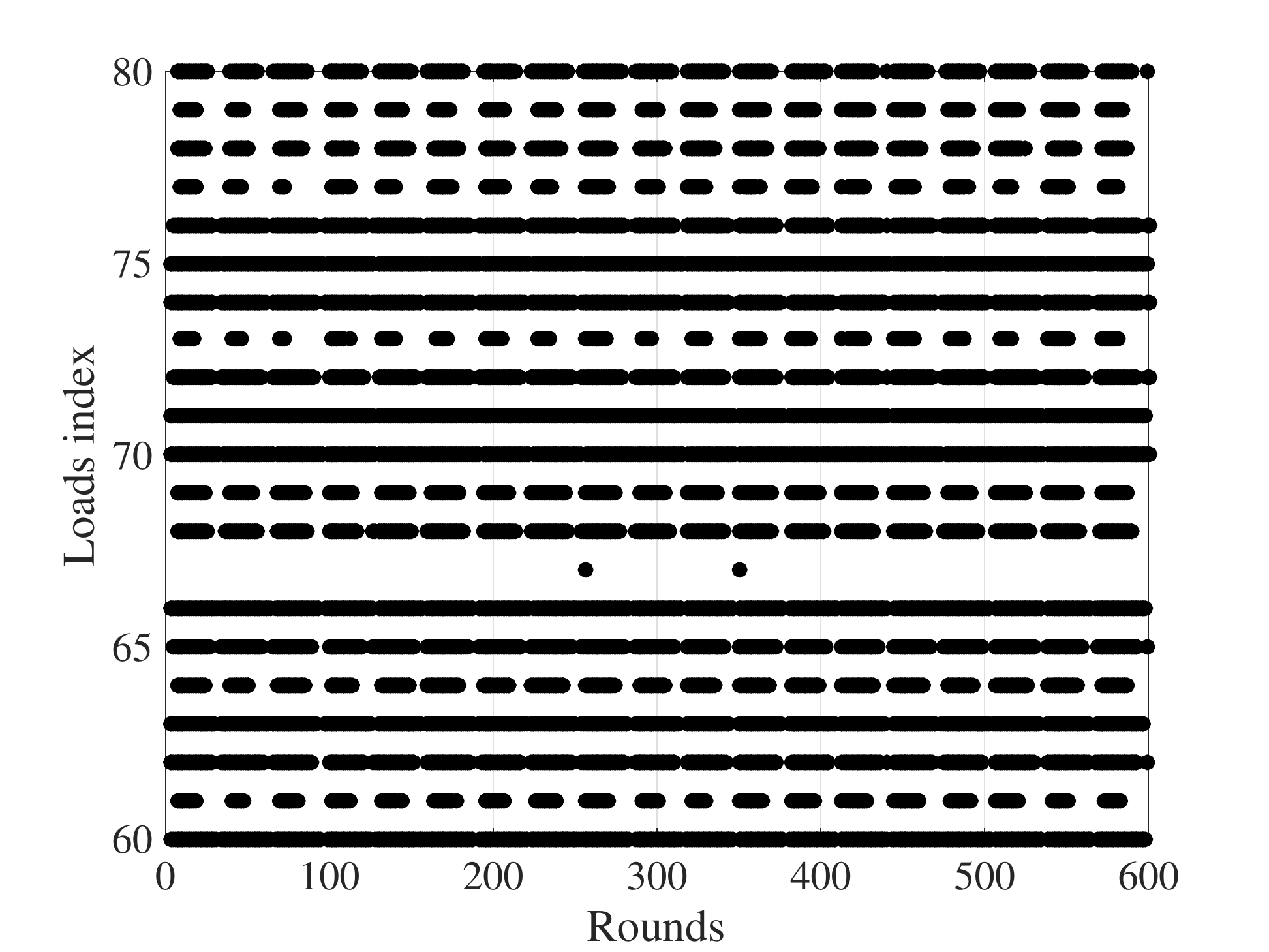}%
\label{fig:sparsity}}

\subfloat[Mean signal value with and without regularization]{\includegraphics[width=0.75\columnwidth,height=0.75\textheight,keepaspectratio]{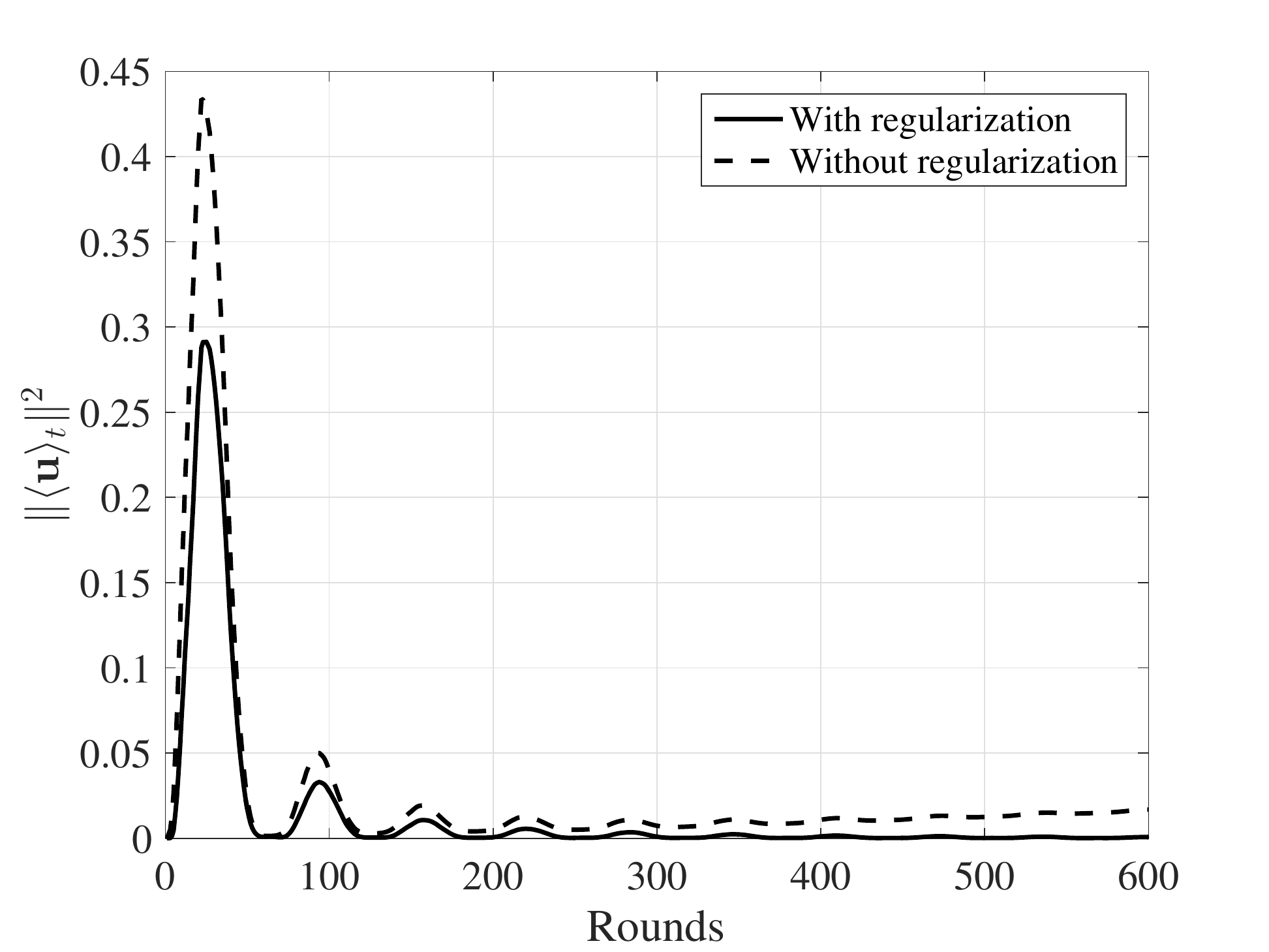}%
\label{fig:mean_reg}}
\caption{Regularization performance of the \texttt{COGD} for setpoint tracking in the full information setting}
\label{fig:fi_reg}
\end{figure}

\begin{figure}[!t]
    \centering
    \includegraphics[width=1\columnwidth]{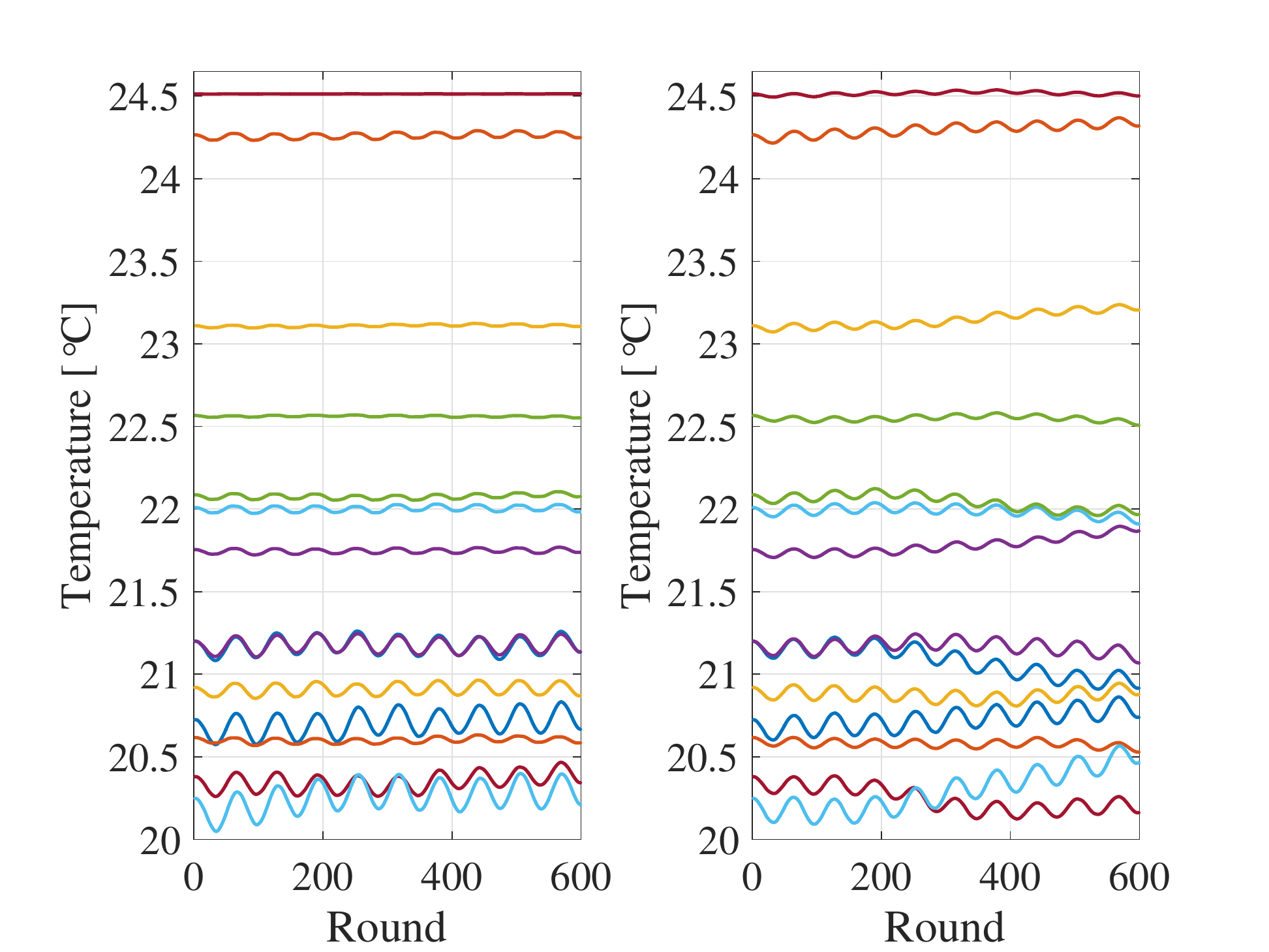}
    \caption{Temperature inside a subset of TCLs in the full information setting. (left) with regularization. (right) Without regularization}
    \label{fig:tcls_temp}
\end{figure}

Figure~\ref{fig:fi_reg} shows the effects of the regularizers in the case of the full information. Note, that these figure shows the combined effect of both regularizers. Figure~\ref{fig:sparsity} shows which loads received a non-zero signal in each round. Figure~\ref{fig:mean_reg} shows that regularizers prevent loads from being dispatched unevenly over time.

Lastly, Figure~\ref{fig:tcls_temp} presents the temperature of a subset of TCLs and shows the combined effect of both regularizers in the full information setting. We observe that for one of the loads the temperature is almost constant reflecting the sparsity of the adjustment signal while for the others, the temperature oscillates in the vicinity of their desired temperature ($\theta_d = \theta(1)$). Without regularization, the temperature drifts away from the desired temperature.

\subsection{Setpoint tracking using EVs storage unit}
\label{sec:evs}
We now apply OCO to setpoint tracking with EVs. To adapt the framework to EVs, some assumptions must be modified. We assume that EVs can either provide or store energy. We also assume that the initial state-of-charge, $S_1$, is at $75\%$ for all EVs. Then, we split the adjustment signal $\mub$ into a charging signal $\mub_c \in [0,1]^N$ where the EV stores energy and hence increases its power consumption, and a corresponding discharging signal $\mub_d \in [-1,0]^N$. We make this assumption to avoid a non-convex loss function. We note that this formulation implies that charging and discharging can occur at the same time. As will be shown later, this unwanted consequence of the convex relaxation is avoided with the sparsity regularizer. 

We let $\cc_{c,t}$ and $\cc_{d,t}$ be the maximum charging and discharging rate at round $t$. The setpoint tracking loss function is written as
\begin{equation*}
\ell_t(\mub_t) = \left(s_t - \cc_{c,t}\Tr \mub_{c,t} - \cc_{d,t}\Tr \mub_{d,t}\right)^2.
\label{f_loss_evs}
\end{equation*}
We modify the mean regularizer for the present non-symmetric formulation because charging and discharging have different effects on the state of charge. The state-of-charge at time $t+1$ for the EVs $i$ is,
\begin{align}
S_{t+1}(i) &= S_t(i) + \frac{h}{B}\left[ \eta^\text{inj}_{s}(i) c_{c,t}(i) \mu_{c,t}(i) \vphantom{\frac{1}{\eta^\text{ext}_{s}(i)}} \right.\nonumber \\
&\quad \left.+ \frac{1}{\eta^\text{ext}_{s}(i)} c_{d,t}(i) \mu_{d,t}(i)  \right],
\label{eq:soc}
\end{align}
where $\eta_s$ the extraction/injection coefficient, $B$ the battery energy capacity and $h$ the length of the time step. We neglect battery leakage in~\eqref{eq:soc} because it is likely to be very small on the faster timescales we are considering. Defining $\langle \mub^w \rangle_t$ as the weighted average of charging and discharging signals for EV, we have
\begin{align}
\langle \mu^w(i) \rangle_t &= \frac{1}{t} \sum_{n=1}^t \left(\eta^\text{inj}_{s}(i) c_{c,n}(i) \mu_{c,n}(i) \vphantom{\frac{1}{\eta^\text{ext}_{s}(i)}} \right. \nonumber\\
&\qquad \left.+ \frac{1}{\eta^\text{ext}_{s}(i)} c_{d,n}(i) \mu_{d,n}(i) \right), \nonumber
\end{align}
\begin{align}
\phantom{\langle \mu^w(i) \rangle_t} &= \frac{1}{t} \left( (t-1)\langle \mu^w_{t-1}(i) \rangle + \eta^\text{inj}_{s} c_{c,t}(i) \mu_{c,t}(i) \vphantom{\frac{1}{\eta^\text{ext}_{s}}} \right.\nonumber\\
&\qquad + \left. \frac{1}{\eta^\text{ext}_{s}} c_{d,t}(i) \mu_{d,t}(i) \right) \label{eq:impact_evs}
\end{align}
for $i=1, \ldots, N$. The mean regularizer is the squared $\ell_2$-norm of the vector $\langle \mub^w \rangle_t$, and prevents the EVs from being dispatched far from their initial state of charge.

For this application example, we only use the full information setting. We do not apply limited feedback because the randomly perturbed update does not ensure that only one of decision variables, $\mub_c$ or $\mub_d$, be non-zero in each time period. 

We set $N=100$ and the time step $h = 1$ minute. The load response is $c_{c,t}(i) = T_\text{charging} + w_t$ and $c_{d,t}(i) = T_\text{discharging} + w_t$ where $w_t\sim \mathrm{N}_{[-1.5,1.5]}(0,0.1)$, a truncated Gaussian noise, models the uncertainty. The setpoint to track is $s_t = 25 \sin (0.1 t)$. For the storage unit, we use $\eta^\text{inj}_{s}(i)= \eta^\text{ext}_{s}(i)=0.85$ for all $i$, a battery capacity of $B=10$ kWh and charging and discharging rates given by $T_\text{charging} =3$~kW and $T_\text{discharging} = 1.5$~kW. These parameters are chosen according to~\cite{ma2013decentralized,tan2014optimal,qin2012optimal}. Finally, we set $\chi = 35$, $\lambda = 46$ and $\rho = 100$. Using our proposed approach, we observe a decrease of $67.08\%$ in the total setpoint tracking loss. Using the regularizers, in average per round, the sparsity is improved by $84.41\%$ and the impact on comfort decreases by $59.71\%$. This results means that the drift in state-of-charge level is significantly reduced.

\begin{figure}[!t]
    \centering
    \includegraphics[width=1\columnwidth]{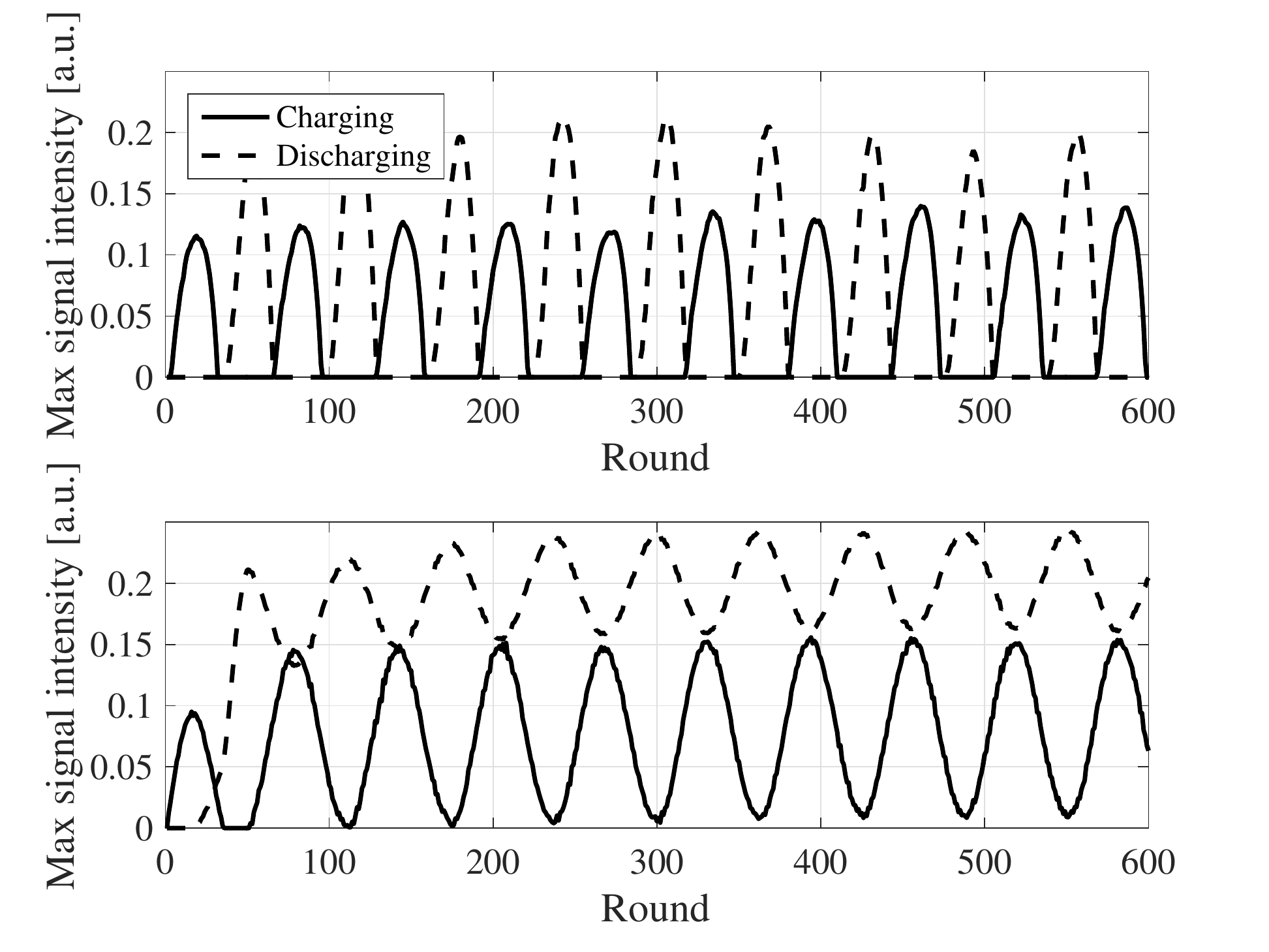}
    \caption{Charging $\mub_c$ and discharging $\mub_d$ signal comparison. (top) With regularization. (bottom) Without regularization.}
    \label{fig:evs_comp}
\end{figure}

Finally, we observe that almost no simultaneous charging and discharging occurs when regularizers are used in the full information setting. Figure~\ref{fig:evs_comp} presents $\max_i |\mu_{c,t}(i)|$ and $\max_i |\mu_{d,t}(i)|$, the maximum signal with and without regularization as a function of $t$. Note that this does not represent the energy in and out of the storage units but only the instruction sent to them. In the top subfigure of Figure~\ref{fig:evs_comp}, as expected when no regularization is used, charging and discharging are dispatched at the same time. Figure~\ref{fig:evs_comp}, bottom subfigure, shows that with the regularizers, simultaneous no charging and discharging almost never occurs. 

\section{Conclusion}
In this work, we propose several online convex optimization algorithms for setpoint tracking using demand response when the load's response is uncertain. We consider three limited feedback extensions: bandit feedback, partial bandit feedback, and Bernoulli feedback. We add two regularizers to our formulation to minimize the number of loads dispatched in each round and to minimize the impact on the loads. We show that each extension's regret has a sublinear upper bound. We apply our model to thermostatically controlled loads and electric vehicles. In our numerical tests, we showed that the limited feedback implementations can still offer good performance while requiring less communication.


\appendix

\subsection{Proof of Theorem~\ref{thm:partial}}
\label{app:partial}

Assume that the aggregate response of the unobserved loads, $\beta_t$, is available in all rounds $t$. Note that all decision variables are known under limited feedback and hence do not change the regularizer, $r(\mub_t)$. 

Similar to $\beta_t$, we define $i_t$ as the effect due to the full information decision variable at round $t$. Then, the loss function $f_t(\mub_t)$ can be rewritten as a function of bandit or full information feedback decision variable only,
\begin{equation}
f_t(\mub_t) = \tilde f_t^F (\beta_t, \mub_t^F) = \tilde f_t^B (i_t, \mub_t^B). \label{eq:loss_f_bernoulli}
\end{equation}
where both $\tilde f_t^F$ and $\tilde f_t^B$ can be computed at round $t$ and each of them are a function of only one type of feedback. Hence, the loss function can be rewritten as,
\begin{equation*}
f_t(\mub_t) = \alpha \tilde f_t^F (\beta_t, \mub_t^F) + (1-\alpha) \tilde f_t^B (i_t, \mub_t^B),
\end{equation*}
for $0<\alpha<1$. Note that $\alpha$ cannot equal zero or one because then the expression would not contain the full $\mub$ vector, and hence would not contain the full regularizer. We do not allow $\alpha$ to depend on $T$. For convenience, we set $\alpha=1/2$. The regret of the partial bandit feedback algorithm, \texttt{PBCOGD}, becomes
\begin{align}
R_T(\texttt{PBCOGD}) &= \sum_{t=1}^T F_t(\mub_t) - F_t(\mub^\ast) \nonumber\\
&= \sum_{t=1}^T \left( f_t(\mub_t) + r(\mub_t^B) + r(\mub_t^F) \right) \nonumber\\
&\quad - \left( f_t(\mub^\ast) + r(\mub^{B\ast}) + r(\mub^{F\ast}) \right) \nonumber\\
&= \sum_{t=1}^T  \frac{1}{2} \left( \tilde f_t^F (\beta_t, \mub_t^F) + \tilde f_t^B (i_t, \mub_t^B) \right) \nonumber\\
&\quad - \frac{1}{2} \left( \tilde f_t^F (\beta^\ast, \mub^{F\ast}) + \tilde f_t^B (i_t, \mub^{B\ast}) \right) \nonumber\\
&\quad+  r(\mub_t^B) + r(\mub_t^F) - r(\mub^{B\ast}) - r(\mub^{F\ast})  \nonumber\\
&= \frac{1}{2}\sum_{t=1}^T \left[ \tilde f_t^F (\beta_t, \mub_t^F) + 2 r(\mub_t^F) \right. \nonumber\\
&\qquad \quad\left. - \tilde f_t^F (\beta^\ast, \mub^{F\ast}) - 2r(\mub^{F\ast})\right] \nonumber\\
& \quad + \frac{1}{2} \sum_{t=1}^T \left[ \tilde f_t^B (i_t, \mub_t^B) + 2 r(\mub_t^B) \right. \nonumber\\
&\qquad \quad \left. - \tilde f_t^B (i_t, \mub^{B\ast}) - 2 r(\mub^{B\ast}) \label{eq:pre_regret_pbcogd}\right] 
\end{align}
The \texttt{PBCOGD} algorithm uses the \texttt{COGD} update on $\tilde f_t^F (\beta_t, \mub_t^F)$ to compute $\mub_{t+1}^F$ and the \texttt{BCOGD} update on $\tilde f_t^B (i_t, \mub_t^B)$ to compute $\mub_{t+1}^B$. Note that, the \texttt{COGD} and \texttt{BCOGD} algorithms are applied to the loss functions in~\eqref{eq:loss_f_bernoulli}. $\beta_t$ and $i_t$ are now seen as new sources of uncertainty by the algorithms. For example, in the setpoint tracking loss case, the full information algorithm tracks a signal given by $s_t + \beta_t$ rather than $s_t$. 

For completeness, the point-wise gradient $\gt^p$ used by the \texttt{BCOGD} update is,
\begin{align*}
\gt^p &= \frac{(N-n)}{\delta} \tilde f_t^B (i_t, \mub_t^B + \delta \vb^p_t) \ \mathbf{v}^p_t \\
&=\frac{(N-n)}{\delta} f_t\left(\mub_t + \delta \begin{pmatrix} \vb^p_t \\ \mathbf{0} \end{pmatrix}\right) \mathbf{v}^p_t,
\end{align*}
as given on line 7 of the algorithm in Figure~\ref{alg:partial}, where $\vb^p_t \sim \Sone^{N-n}$ for all $t$.

Finally, substituting in the regret bounds for \texttt{COGD} and \texttt{BCOGD}, we can rewrite~\eqref{eq:pre_regret_pbcogd} as
\begin{align*}
R_T(\texttt{PBCOGD}) & = \frac{1}{2} R_T(\texttt{COGD}; \tilde f_t^F) + \frac{1}{2} R_T(\texttt{BCOGD}; \tilde f_t^B)\\
& \propto O\left(T^{\frac{3}{4}}\right). \hskip0.52\columnwidth \text{\QED}
\end{align*}

\subsection{Proof of Theorem~\ref{thm:bernoulli}}
\label{app:bern}
Let $T_{B}$ be number of times the player receives bandit feedback over $T$ trials. Then, $T_{B} \sim \mathrm{Binomial}(T,p)$ and the number of times the player receives full feedback is $T-T_B$. 
Using the approach of the proof of Lemma~\ref{lem:tun-comid}~\cite{duchi2010composite,oco_irep2017} with $\mathcal{R}(\cdot) = \frac{1}{2}\norm{\cdot}_2^2$ and of the proof of Theorem~\ref{thm:bcogd}~\cite{oco_irep2017}, we observe that
\begin{align*}
\E[R_T(\texttt{BerCOGD})] \leq& \E\left[r(\mub_1) + \frac{D^2}{2}\left(\frac{1}{\eta_F} + \frac{1}{\eta_B} \right)\right. \\
&\qquad \left. + \frac{1}{2} \sum_{t=1}^T w_t + \sum_{t=1}^T z_t \right]
\end{align*}
where,
\begin{align*}
w_t &= \begin{cases}
\eta_F G^2, \text{ if } I_t = 0 \\
\frac{\eta_B N^2}{\delta^2} B^2, \text{ if } I_t = 1
\end{cases} \\
z_t &=\begin{cases}
0, \text{ if } I_t = 0 \\
2\delta DL + 2\delta L, \text{ if } I_t = 1.
\end{cases}
\end{align*}
We remind the reader that $I_t=1$ in the case of bandit feedback and $I_t=0$ in the case of full feedback. Therefore, by linearity of $R_T$ and taking into consideration the initialization steps, we re-express the conditional expected regret of the algorithm as
\begin{align*}
\E[R_T(\texttt{BerCOGD})] &= \E\left[ \E \left[ \left. R_T \right| T_{B} + 1\right] \right] \\
&= \E \left[ \E\left[  \left. R_T(\texttt{BCOGD}) \right| T_B + 1\right]\right]\\
&\qquad + \E\left[ \E \left[  \left. R_T(\texttt{COGD}) \right| (T-T_{B} +1)\right]\right] \\
\end{align*}
since $r(\mub_1) = 0$. Then, using Lemma~\ref{lem:tun-comid} and Theorem~\ref{thm:bcogd} for the inner conditional expectations, we have,
\begin{equation*}
\E[R_T(\texttt{BerCOGD})] \leq \E\left[ A_1 \left(T-T_B + 1\right)^{\frac{1}{2}} + A_2 \left(T_{B}+1\right)^{\frac{3}{4}} \right]
\end{equation*}
with $A_1 = DG\chi_1$ and $A_2 = DBN\chi_2 + 2DL + 2L$ and where we used
\begin{align*}
\eta_F &=  \frac{D \chi_F}{G \left(T-T_B + 1\right)^\frac{1}{2}}, \\
\eta_B &= \frac{D \chi_B}{B N (T_B+1)^{\frac{3}{4}}}, \\
\delta_B &= \frac{1}{(T_B+1)^{\frac{1}{4}}},
\end{align*}
to bound the regret in the inner expected value.

$(T_B +1)^{\frac{3}{4}}$ and $\left(T-T_B + 1\right)^\frac{1}{2}$ are concave functions for $0 \leq T_{B} \leq T$. Hence, by Jensen's inequality, we have
\begin{align}
\E[R_T(\texttt{BerCOGD})] &\leq A_1 \E[(T-T_B + 1)]^{\frac{1}{2}} + A_2 \E[T_B + 1]^{\frac{3}{4}} \nonumber \\
&\leq A_1 (T - Tp + 1)^{\frac{1}{2}} + A_2 (Tp + 1)^{\frac{3}{4}} \nonumber \\
&\leq A_1 T^{\frac{1}{2}} + A_2 T^{\frac{3}{4}} p^{\frac{3}{4}} + A_1 + A_2 \nonumber
\end{align}
since $T_B$ is a binomial random variable and $0 \leq p \leq 1$. Setting $p \leq a/T^{\frac{1}{3}}$ leads to the given expected regret upper bound. Lastly, $\eta_F$, $\eta_B$ and $\delta_B$ are set according to the number of feedback rounds $T_B$ sampled in the initialization step. Note that the initialization rounds are necessary to ensure that $\eta_F$, $\eta_B$ and $\delta_B$ are defined. \QED


\ifCLASSOPTIONcaptionsoff
  \newpage
\fi



%







\end{document}